\documentclass[10pt]{amsart}
\usepackage[cp1251]{inputenc}
\usepackage[english,russian]{babel}
\usepackage{amsmath}
\usepackage{amssymb}
\usepackage{amsfonts}
\usepackage{wrapfig}
\usepackage{graphicx,color}

\newtheorem*{theorem}{Theorem}

\begin{document}
	\renewcommand{\refname}{References}
	
	\thispagestyle{empty}
	
	\title[On periodic groups isospectral to $A_7$]{On periodic groups isospectral to $A_7$}
	\author{{A. MAMONTOV}}%
	
	\thanks{\sc Mamontov, A.,
		On periodic groups isospectral to $A_7$}
	\thanks{\copyright \ 2018 Mamontov A}
	\thanks{\rm This research is supported by RSF (project 14-21-00065)}
	\maketitle {\small
		\begin{quote}
			\noindent{\sc Abstract. } The spectrum of a periodic group $G$ is the set $\omega(G)$ of its element orders. Consider a group $G$ such that $\omega(G)=\omega(A_7)$. Assume that $G$ has a subgroup $H$ isomorphic to $A_4$, whose involutions are squares of elements of order $4$. We prove that either $O_2(H) \subseteq O_2(G)$ or $G$ has a finite nonabelian simple subgroup.\medskip
			
		\end{quote}
	}

\section{Introduction}

Let $\mathfrak{M}$ be a set of periodic groups and $G \in \mathfrak{M}$.  The {\it spectrum} of $G$ is the set $\omega (G)$ of its elements orders. Two groups with the same spectrum are called {\it isospectral}. We say that $G$ is {\it recognizable by spectrum in} $\mathfrak{M}$, if any group from $\mathfrak{M}$ is isomorphic to $G$, whenever it is isospectral to $G$.

Many finite simple groups are known to be recognizable by spectrum in the class of finite groups \cite{fin_review}, in particular, $A_7$ \cite{OCn}. It is clear that a group $G$ is recognizable by spectrum in the class of periodic groups, if it is recognizable among finite groups, and any group isospectral to $G$ is locally finite. Verification of the last condition is related to the 
Burnside problem \cite{lmreviewe}.

In the paper we prove

\begin{theorem} Assume $\omega(G)=\omega(A_7)$ and $G$ has a subgroup $H$ isomorphic to $A_4$, whose involutions are squares of elements of order  $4$. Then either  $O_2(H) \subseteq O_2(G)$ or $G$ has a finite nonabelian simple subgroup.
\end{theorem}

Existence of such simple subgroup $K \leq G$ allows to investigate the structure of $G$ using the inclusion $C_K(i) \leq C_G(i)$ (where $i$ is an involution from $K$) and Shunkov's theorem that a periodic group with a finite centralizer of involution is locally finite \cite{shu1972e}.

We conjecture that $G \simeq A_7$.


\section{Notations and preliminary results}

Speaking of computations we refer to computations in GAP \cite{gap} using coset enumeration algorithm. The following notations are used: $\Gamma_n=\Gamma_n(G)$ is the set of elements of order $n$ from 
$G$; $\Delta=\Delta(G)=\{x^2\mid x \in \Gamma_4 \}$, $O_2(G)$ is the largest normal 2-subgroup of $G$. Write $A:B$ for an extension of $A$ by $B$, and $C_n$ for a cyclic group of order $n$. The center of $G$ is denoted as $Z(G)$. Writing $a \sim b$ for $a,b \in G$ we mean that orders of $a$ and $b$ are equal.

Groups in the work are mostly presented by generators and relations. Any word identity in a group also holds in its homomorphic images. So it is convenient to keep same notations for generators of both groups.
	
For visualization some relations are drawn in the form of labeled graph. Vertex label represents some group element, and the corresponding number is its order. Edge label shows relations between the corresponding vertices (number stands for the order of the vertex element product), or subgroup, which they generate. An edge between commuting vertices is drawn as dashed line.

In this section let $G$ be a group whose element orders are not greater than $7$. The goal of this section is to describe small subgroups, which may be in $G$. For further discussion, we distinguish some elements, which are written as words of generators, and establish
some relations between them.

In the following lemma we describe $(2,3)$-generated subgroups of $G$.

\begin{wrapfigure}[120]{l}[340pt]{20mm}
	\includegraphics[width=20mm]{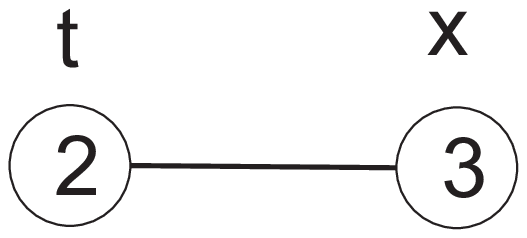}
\end{wrapfigure}

{\bf Lemma 2.1} {\it Consider an involution $t$ and an element $x$ of order $3$ in $G$. Let $K=\langle t,x \rangle$. Then one of the following holds:

\medskip

\begin{itemize}

\item[(a)] $(xt)^{2}=1$ and $K \simeq S_3$;

\item[(b)] $(xt)^{3}=1$ and $K \simeq A_{4}$;

\item[(c)] $(xt)^{4}=1$ and $K \simeq S_{4}$;

\item[(d)] $(xt)^{5}=1$ and $K \simeq A_{5}$;

\item[(e)] $(xt)^{6}=1$, and $K$ is a homomorphic image of $(C_k \times C_k) \leftthreetimes C_6$, where $k$  is the order of $[x,y]$;

\item[(f)] $(xt)^{7}=[x,t]^4=1$ and $K \simeq L_2(7)$.

 \end{itemize}
}

\begin{proof}[Proof]

Items (a)-(d) are well-known.  

(e) Let $(xt)^6=1$ and $k$ be the order of $[x,t]$. Set
$a=[x,t]=x^{-1}txt$. Then $a^t=tx^{-1}txt\cdot t=tx^{-1}tx=a^{-1}$.
Further, $aa^xa^{x^2}=(ax^{-1})^3=x^{-1}(tx)^6x=1$ and
$a^{xt}=a^{txa}=(a^{-1})^{xa}$, therefore $\langle a,a^x \rangle
\unlhd \langle x,t \rangle$. Direct computations show that 
$[a,a^x]=(tx)^6=1$.

(f)  Let $(xt)^7=1$. Then $K$ is a homomorphic image of 
$K(j)=\langle x,t | 1=x^3=t^2=(xt)^7=[x,t]^j \rangle$, where $j
\in \{4,5,6,7 \}$. Computations show that $K(5) \simeq 1$;
$K(6) \simeq K(7) \simeq L_2(13)$ and therefore has an element of order $13$, which is not possible. Hence, $K \simeq K(4) \simeq
L_2(7)$.

The lemma is proved.

\end{proof}

Further $F_n$ will denote the Frobenius group $\langle x,t| 1=(xt)^6=[x,t]^k \rangle$ from item (e), where $k \in \{ 4,5,6,7 \}$; and $n=6k^2$ is the order of the group. In the following Lemma 2.3 we list some properties of $F_n$ that will be used later.

\begin{wrapfigure}[120]{l}[340pt]{20mm}
	\includegraphics[width=20mm]{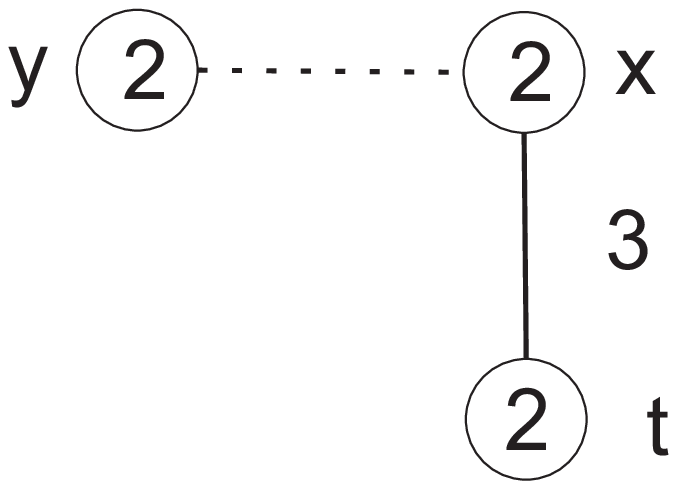}
\end{wrapfigure}

{\bf Lemma 2.2} {\it 
	Let $x,y,t \in \Gamma_2$, $[x,y]=1$ and $\langle x,t \rangle
	\simeq S_3$. Then the group $\langle x,y,t \rangle$ is finite and $yt \not
	\in \Gamma_7$. Moreover, if $yt \in \Gamma_5$, then $\langle x,y,t
	\rangle \simeq A_5$.}

\begin{proof}[Proof]

Note that $\langle x,y,t \rangle$ is a homomorphic image of $\langle x,y,t \mid
1=x^2=y^2=t^2=(xy)^2=(xt)^3=(yt)^i=(xyt)^j \rangle$, where $i,j
\in \{4,5,6,7 \}$ are selected in correspondence with Lemma 2.1. Now the statement is easily verified by computations.

The lemma is proved. \end{proof}

{\bf Lemma 2.3} {\it Let $z=t^{xt}, b=z^xz,
	u=z^{x^2z}=t^{[x^{-1},t][x,t]}$ be the corresponding elements of $F_{294}$. Then $b^x=b^{4}$, $F_{42} = \langle
	x,z | 1=x^3=z^2=(xz)^6=b^7, b^x=b^4 \rangle$ is the Frobenius group of order $42$, and $[u,x]=1$.}

\begin{proof}[Proof]

Let $a=t^xt$ as in the proof of Lemma 2.1. Then $z=ta$ and
$$b=x^{-1}taxta=aa^{xt}a=a(a^{-1})^{xa}a=(a^{-1})^xa^2,$$

$$b^4=(a^3)^xa=x^{-1} \cdot (x^{-1}txt)^3 \cdot x \cdot x^{-1}txt
=$$

$$=xtxtx^{-1}txtx^{-1}tx^{-1}t,$$

$$b^x=(a^2(a^{-1})^x)^x=x^{-1} \cdot x^{-1}txtx^{-1}txt \cdot x^{-1}
\cdot tx^{-1}tx \cdot x^2=xtxtx^{-1}txtx^{-1}tx^{-1}t.$$

Note that a map 

$$x \rightarrow
(1,2,3)(4,5,6)(7,8,9)(10,11,12)(13,14,15)(16,17,18)(19,20,21)$$
$$t \rightarrow (2,4)(3,5)(6,7)(8,10)(9,11)(12,13)(14,16)(15,17)(18,19)$$
can be extended to an isomorphism between $\langle x,t \rangle $ and the corresponding permutation group.

Now the rest of the proof is straightforward. \end{proof}

This notation $F_{42}$ for the corresponding Frobenius group will be used further. An important property of subgroups of $G$, which are isomorphic to $F_{42}$, is described in the following section. 


\section{Subgroups isomorphic to $D_{10}$ and $F_{42}$}

	
The goal of the section is to prove two statements. Writing $F_{42}=\langle z,x \rangle$ we assume that elements $z$ and $x$ satisfy the relations, which were used to define $F_{42}$ in Lemma 2.3.

\begin{wrapfigure}[120]{l}[340pt]{20mm}
	\includegraphics[width=20mm]{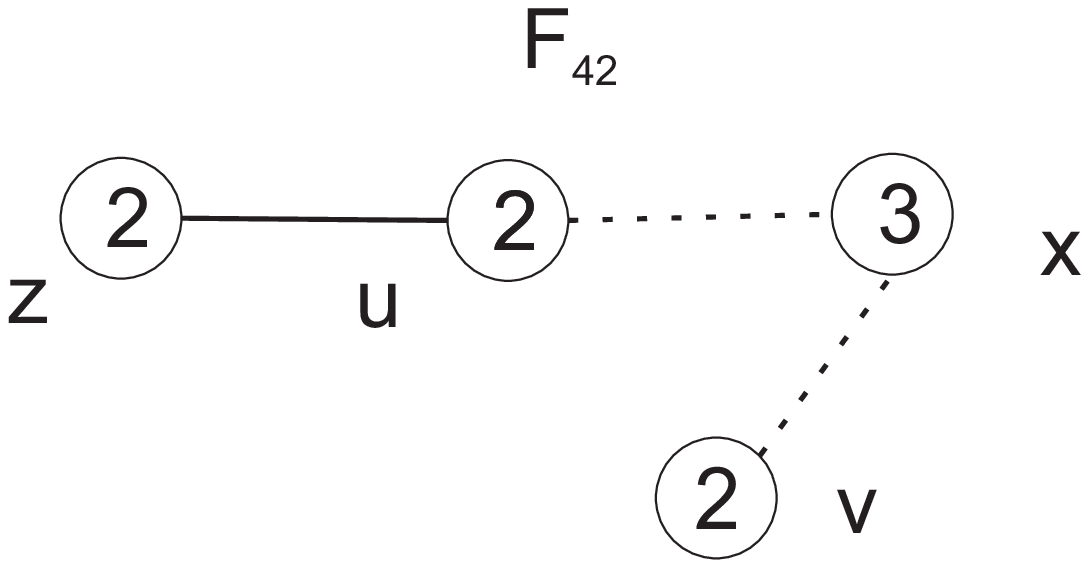}
\end{wrapfigure}

{\bf Statement 1} {\it Assume that $G$ has no elements of order greater than $7$. Let $F_{42}= \langle z,x \rangle$ be a subgroup of $G$, $u=z^{x^2z}$, and $v$ be an involution in $G$ such that $[x,v]=1$. Then either $v=u$ or $G$ has a subgroup isomorphic to $L_2(7)$.}

The proof is preceded by two lemmas with the same assumptions. Denote by $B$ the following set of relations, which are assumed, and also shown on the graph:

$$B=\{ v^2, z^2, x^3, (xz)^6, b^7, b^xb^{-4}, [x,v] \}, \text{where } b=z^xz.$$

\begin{wrapfigure}[120]{l}[340pt]{20mm}
	\includegraphics[width=20mm]{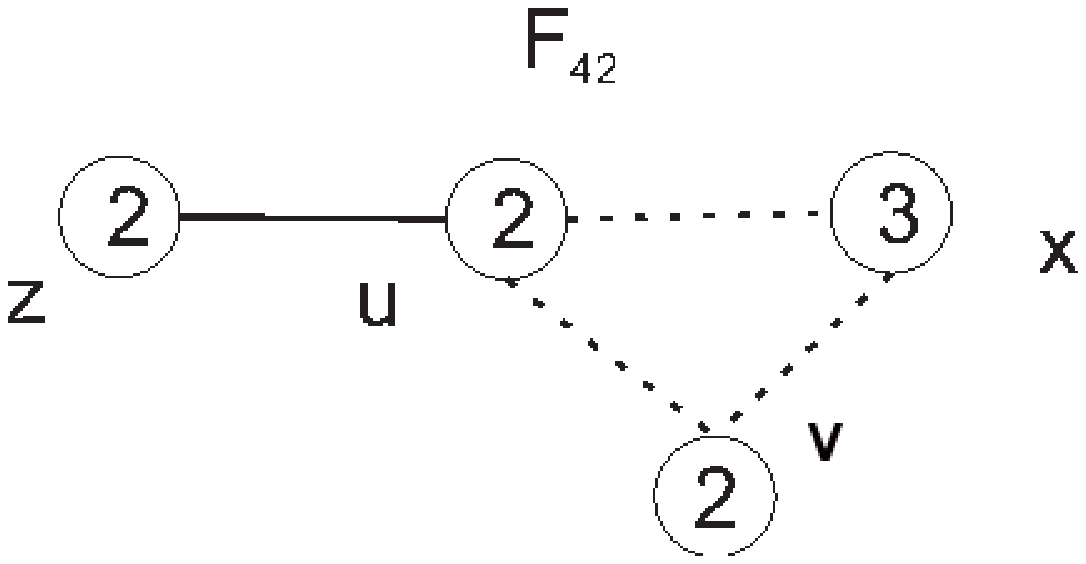}
\end{wrapfigure}

{\bf Lemma 3.1} {\it If $[u,v]=1$, then $v=u$.}

\begin{proof}[Proof]
	
Note that the subgroup $\langle v,z,x \rangle$ is a homomorphic image of $G(i_1,i_2,i_3,i_4) = \langle v,z,x| B \cup \{ [u,v], (zv)^{i_1}, (v^{{x^z}}v)^{i_2}, (vx^z)^{i_3}, (xzv)^{i_4} \} \rangle$, where $i_1,\ldots,i_4 \in \{ 4,5,6,7 \}$. Computations show that either the order of  $G(i_1,i_2,i_3,i_4)$ divides $42$, or $v$ centralizes $F_{42}$,  which is not possible.
The lemma is proved. \end{proof}

\begin{wrapfigure}[120]{l}[340pt]{20mm}
	\includegraphics[width=20mm]{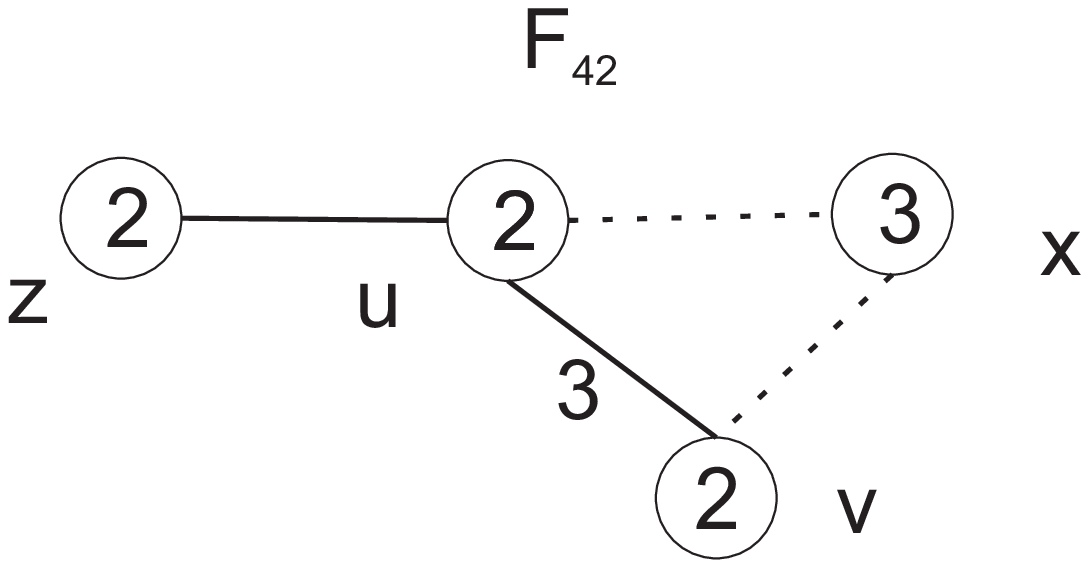}
\end{wrapfigure}

{\bf Lemma 3.2} {\it  If $(uv)^3=1$, then either $v=u$ or $G$ has a subgroup isomorphic to $L_2(7)$.}

\begin{proof}[Proof]
	
Assume the contrary. Then by Lemma 2.1.e the order of $v \cdot x^z$ and $v \cdot xx^z$ is not 7. Therefore $\langle x,z,v \rangle $ is a homomorphic image of $G(j_1,j_2,i_1,i_2,i_3) = \langle x,z,v| B \cup \{ (uv)^3, (vx^z)^{j_1}, (xx^zv)^{j_2}, (vz)^{i_1}, (vzx)^{i_2}, (xx^zxv)^{i_3} \} \rangle$, where $j_1, j_2 \in \{4,5,6\}$ and $i_1,i_2, i_3 \in  \{4,5,6,7\}$. Computations show that indexes  $|G(j_1,j_2,i_1,i_2,i_3): \langle x \rangle|$ are not greater than $14$, which means that $v \in \langle x,z \rangle$. A contradiction. The lemma is proved. \end{proof}

\begin{proof}[Proof of Statement 1] 

Assume the contrary. The dihedral subgroup $\langle u,v \rangle$ centralizes $x$ and so does not contain elements of orders $4$, $5$ and $7$. If the order of $uv$ is even, then let $i$ be the involution from the center of $\langle u,v \rangle$. By Lemma 3.1 $i=u$ and hence $v=u$, a contradiction. If the order of $uv$ is 3, then obtain a contradiction by Lemma 3.2. The statement is proved. \end{proof}

We will also need to following corollary.

\begin{wrapfigure}[120]{l}[340pt]{20mm}
	\includegraphics[width=20mm]{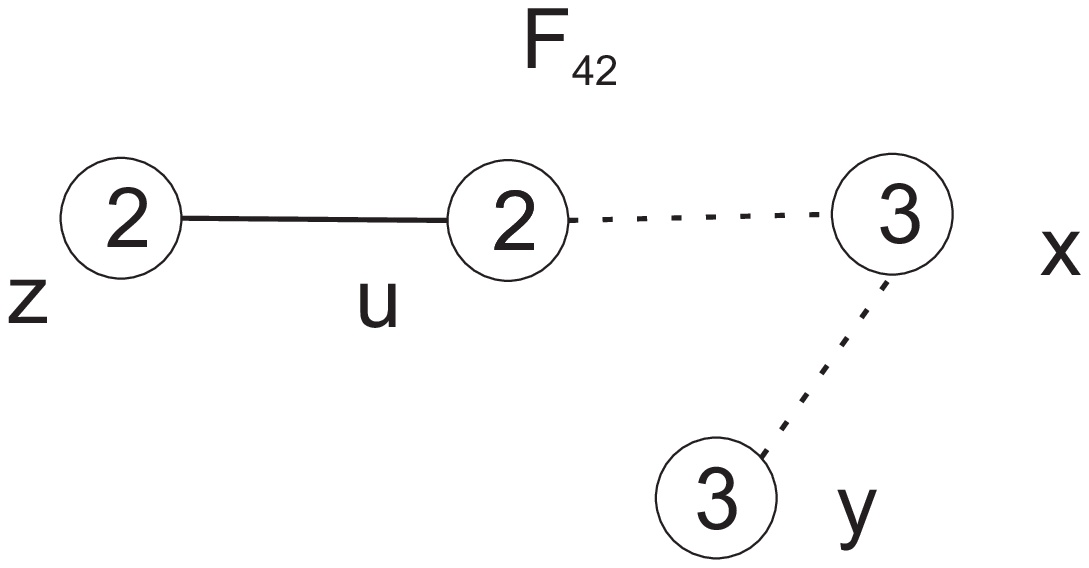}
\end{wrapfigure}

{\bf Lemma 3.3} {\it Let $F_{42}= \langle z,x \rangle$ be a subgroup of $G$, and let $y$ be an element of order $3$ from $C_G(x)$. Then $[y,u]=1$, where $u=z^{x^2z}$.}

\begin{proof}[Proof] There are no elements of order $12$ in $G$, and so there are no elements of order $4$ in $\langle u,y \rangle$. By Statement 1 there are no involutions in $\langle u,y \rangle$ other than $u$. By Lemma 2.1 we obtain $[u,y]=1$. The lemma is proved. \end{proof}

{\bf Statement 2} {\it Assume that $4,5 \in \omega(G)$ and $G$ has no elements of order greater than $7$. Let $a \in \Delta$. If $a$ is inverting an 
element of order $5$ and $C_G(a)$ contains an involution $t \not = a$, then $G$ has a nonabelian finite simple subgroup. }

\begin{proof}[Proof]

Let $a\in \Delta$. Take $b \in \Gamma_4$ such that $b^2=a$. By assumption there is an involution $c$ such that $ac \in \Gamma_5$.

First we use relations, which define $\langle b,c \rangle$, to prove that $G$ has a Frobenius subgroup $K \simeq F_{20}$, which contains $a$. Then 
we prove that $\langle K,t \rangle$ is finite and satisfies the statement conclusion.

A group $\langle b,c \rangle$ is a homomorphic image of order $\geq 20$ of a group $G(i_1,i_2,i_3,i_4,i_5)=\langle b,c | 1=b^4=c^2=x^5=(xb)^{i_1}=[x,b]^{i_2}=(a^xb)^{i_3}=(bxb^3x^2)^{i_4}=((cb^2cb)^2cb^{-1})^{i_5} \rangle $, where $a=b^2$, $x=ac$; and $i_1,\ldots,i_5 \in \{4,5,6,7\}$. Computations show that 
$G(5,4,5,6,6)\simeq A_6$, $G(7,5,7,4,5)\simeq L_3(4)$, which is not possible; $H=G(4,5,4,5,4) \simeq (C_5 \times C_5):C_4$. For other values of the parameters the order of the corresponding group is not greater than $20$.

In both suitable cases: when $\langle b,c \rangle$ is isomorphic to $H$ and when $\langle b,c \rangle$ is its homomorphic image $F_{20}=C_5:C_4$ we may 
choose $d$ such that $d^4=1, b=dad^2$ and $\langle a,d \rangle$ is a subgroup of $G$ isomorphic to $F_{20}$. Indeed, in $H$ consider $d=b^{cbcaca}$.

In further proof we use the assumption that $C_G(a)$ has an involution $t \not =a$. By \cite[Lemma 9]{lmmj2014e} such involution $t$ may be selected in the normalizer of $b$.

Therefore, $\langle d,a,t \rangle$ is a homomorphic image of $G(w,i_1,i_2,i_3,i_4) = \langle d,a,t|1=d^4=a^2=t^2=[a,t]=w=(td)^{i_1}=(dat)^{i_2}=(d^2t)^{i_3}=(d^2at)^{i_4}, a=b^2 \rangle$, where $b=dad^2$, $w \in \{ b^{-1}b^t, bb^t \}$, $i_1,i_2,i_3,i_4 \in \{ 4,5,6,7 \}$. Computations show that either the orders of these groups are not greater than 20, or they are isomorphic to $S_5$ or $S_6$, which is not possible. For example, $G(b^{-1}b^t,5,5,6,6) \simeq S_6$. 

The statement is proved.  \end{proof}


\section{Subgroups isomorphic to $A_4$}

The main goal of this section is to prove the theorem. So we assume that $G$ satisfies the assumptions of Theorem, and $H$ is the corresponding subgroup, which is isomorphic to $A_4$.

Throughout the section,  $a\in \Delta$  and $x$ denote generators of $H$, which satisfy the following relations:  $x^3=a^2=(ax)^3$. Then $O_2(H)=\langle a,a^x \rangle$. Also in this section we
assume that $G$ is a counterexample to the theorem, and in particular it does not contain finite nonabelian simple subgroups. As a consequence, by Lemma 2.1, for every involution $i \in G$ and element $y \in G$ of order $3$ the order of $iy$ divides $4$ or $6$.

Proof of Theorem is based on the following Baer-Suzuki argument for $p=2$, that was obtained in \cite{bs2014e}.

{\bf Lemma 4.1} {\it
	Let $G$ be a group of period $n=4k$, where $k$ is odd, containing an involution $i \in G$. If any two elements in $i^G$ generate a $2-$group, then $\langle i^G \rangle$ is a normal $2-$subgroup in $G$.
}

Therefore to prove Theorem, it is sufficient to show that if $t$ is an involution, then either $(at)^4=1$ or $\langle t,a,x \rangle$ contains a finite nonabelian simple subgroup. This strategy is implemented in the following lemmas, where we consequently show that in corresponding situations our involution $t$ cannot be a counterexample.

Possibilities for the subgroup $\langle t,x \rangle$ are restricted by Lemma 2.1. First of all we consider the case, when  $\langle t,x \rangle$ is a homomorphic image of $F_{294}$, in particular, when $[t,x]=1$.

\begin{wrapfigure}[120]{l}[340pt]{20mm}
	\includegraphics[width=20mm]{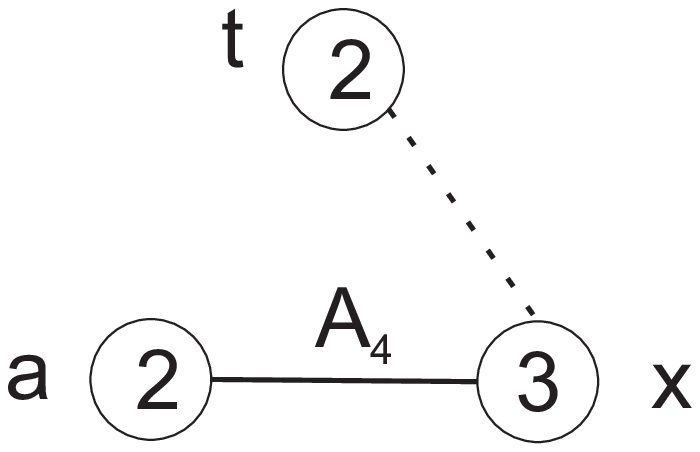}
\end{wrapfigure}

{\bf Lemma 4.2} {\it
	If $t$ is an involution and $[t,x]=1$, then $(at)^4=1$.
}

\begin{proof}[Proof] 

Note that $\langle t,a,x \rangle$ is a homomorphic image of 
$G(i_1,i_2,...,i_9)=\langle a,t,x|
1=a^2=t^2=x^3=(ax)^3=[t,x]=(ta)^{i_1}=(ta^x)^{i_2}=(ata^x)^{i_3}=(a^x(at)^2)^{i_4}=[a^x,(at)^2]^{i_5}=(a^x(at)^3)^{i_6}=(a(ta^x)^3)^{i_7}=(atx)^{i_8}=((at)^2x)^{i_9}
\rangle$, where $i_1,...,i_9 \in \{ 4,5,6,7 \}$. Computations show that this group is finite and is either isomorphic to $S_5$, or is a $\{2,3
\}$-group, and in the last case $O_2(\langle a,x \rangle) \subseteq O_2( \langle a,t,x \rangle)$.

The lemma is proved. \end{proof}

Note that relations that we use are taken from \cite[Lemma 7.A]{lmmj2014e}.

{\bf Lemma 4.3} {\it For any involution $t \in G$ the order of $tt^x$ is not $7$.}

\begin{proof}[Proof] 

Assume the contrary. 
By Lemmas 2.1 and 2.3 there is an involution $z$ such that $\langle z,x \rangle$=$F_{42}$. By Statement~1 $C_G(x)$ has only one involution $u=z^{x^2z}$. By Lemma~4.2 we have $(ua)^4=(ua^x)^4=1$. 

Therefore $\langle u,a,x \rangle$ is a homomorphic image of $L_0=\langle u,a,x |
1=u^2=a^2=x^3=(ax)^3=[u,x]=(ua)^4=(ua^x)^4 \rangle$. Computations show that this group is finite:
$|L_0|=192$ and $z_0=(x^2uaxua)^2$ is the central involution of $L_0$. In particular, $z_0 \in C_G(x)$.

If $z_0=u$, then looking at $\overline{L_0}=L_0/ \langle uz_0 \rangle \simeq A_4$ we see that $\overline{u}=\overline{1}$, which is not possible.

Therefore by Lemma 3.1 we may further assume that $z_0=1$. Hence $\langle u,a,x \rangle$ is a homomorphic image of $L_1=L_0/ \langle z_0=1 \rangle$. In $L_1$ let $i=(ua)^2$. Then $i^2=(ix)^3=[u,i]=1$.

First assume that $i \not = 1$. Then $\langle i,x \rangle \simeq A_4$, and we substitute $a$ by $i$ and consider the group $\langle i,z,x \rangle$
to simplify the relations. Observe that the following relations hold: $R= \{ i^2, z^2, x^3, (xz)^6, b^7, b^xb^{-4}, (ix)^3, (ui)^2 \}$, where $b=z^xz$, $u=z^{x^2z}$. 

Let $y=z^i$. By Lemma~2.1 for $\langle y ,x \rangle$ one of the following possibilities holds:

\begin{itemize}
	
	\item $(yx)^6=(y^xy)^7=1$. 
	
	Let $v=y^{[x^{-1},y][x,y]}$. By Lemma 2.3 $[v,x]=1$ and by Statement~1 $u=v$. Computations show that the order of  $G(k_1)=\langle i,z,x | R \cup \{ (iz)^{k_1}, (yx)^6, (y^xy)^7, uv^{-1} \} \rangle$ for $k_1 \in \{5,7 \}$ divides $3$. Note that $G(4)=G(6)=G(2) \simeq C_7:(C_2 \times A_4)$ has an element of order $14$, and the corresponding homomorphic image without it is just $\langle x,z \rangle$. A contradiction.
	
	\item $(yx)^6=(y^xy)^5=1$. 
	
	Let $v=y^{xy(x^{-1}y)^2}$. Then $[v,x]=1$ and by Statement~1 $u=v$. Computation show that the order of $G(k_1)=\langle i,z,x | R \cup \{ (iz)^{k_1}, (yx)^6, (y^xy)^5, uv^{-1} \} \rangle $ for $k_1 \in \{4,5,6,7 \}$ divides $24$. A contradiction.
	
	\item $(yx)^6=(y^xy)^6=1$. 
	
	Elements $h=(y^xy)^2$ and $k=[h,x]$ of the group $\langle y,x | y^2=x^3=(yx)^6=(y^xy)^6=1 \rangle$ satisfy the following conditions: $\langle h,x \rangle \simeq 3^{1+2}$ is the extraspecial group of exponent 3 and order $27$, and $k$ is the element, generating its center. Moreover, the involution $y$ inverts both $h$ and $k$. 
	
	Note that $uy=uizi\sim iuiz=uz \in \Gamma_7$. By Lemma 3.3 $[u,k]=1$. Therefore the element $(uy)^2$ of order $7$ centralizes $k$. It follows that $k=1$. Changing $k$ to $h$, in an analogous way we obtain that $h=1$. Therefore this case is a particular case of the following one.
	
	\item $(yx)^6=(y^xy)^4=1$. 
		
	Let $v=yx^{-1}yxyxyx^{-1}y$. Then $[v,x]=1$ and by Statement~1 either $u=v$, or $v=1$. Computations show that the order of $$G(k_1,\epsilon)=\langle i,z,x | R  \cup \{ (iz)^{k_1}, (yx)^6, (y^xy)^4, u^{\epsilon}v \} \rangle$$ for $k_1 \in \{4,5,6,7 \}, \epsilon \in \{0,1 \}$ divides $24$. A contradiction.
	
	\item $(yx)^4=1$. 
	Let $v=y^{xy}$. Then $x^v=x^{-1}$. Therefore $vuv \in C_G(x)$ and by Statement~1 $u=vuv$. Computations show that the order of  $$G(k_1,k_2)=\langle i,z,x | R \cup \{ (iz)^{k_1}, (izx)^{k_2}, (yx)^4, (uv)^2 \} \rangle$$  divides $2$ for $k_1, k_2 \in \{4,5,6,7 \}$. A contradiction.

\end{itemize}

In this way we obtain that $i=1$. Then $(ax)^3=(ua)^2=1$, changing $i$ to $a$ and repeating arguments above we obtain that $a=1$. A contradiction. The lemma is proved. \end{proof}

\begin{wrapfigure}[120]{l}[340pt]{20mm}
	\includegraphics[width=20mm]{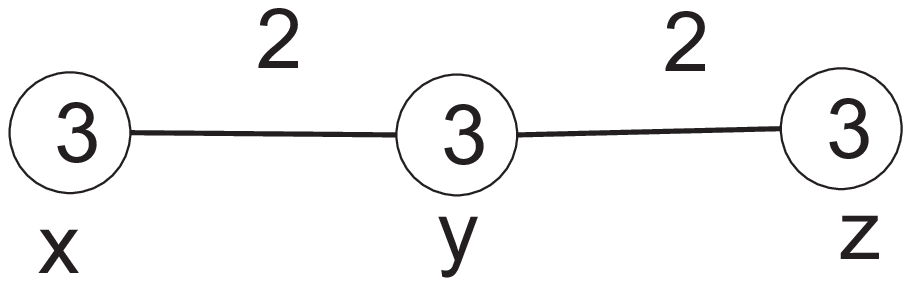}
\end{wrapfigure}

{\bf Lemma 4.4} {\it Let $x,y,z$ be elements of order $3$ such that $(xy)^2=(yz)^2=1$ and $H=\langle x,y,z \rangle$. Then $H$ is an extension of a $2$-group by a group isomorphic to $A_4$ and the following relations hold: $R=\{(xz)^6, (x^{-1}z)^4, (xyz)^6, (xyzxz)^4, (xzy)^6, (x^{-1}z^y)^4 \}$ and $(y^{-1}zx)^4=1$.}

\begin{proof}[Proof] Note that  $x \cdot zy$ and $xzy= xy \cdot z^y$ are products of an involution and an element of order $3$, hence their orders divide $4$ or $6$. 
Let's use the following notations for sets of relations in the alphabet $ \{x,y,z \}$.

$$B= \{ x^3, y^3, z^3, (xy)^2, (yz)^2 \}, $$

$$R(k_1,\dots,k_6)= \{ (xz)^{k_1}, (x^{-1}z)^{k_2}, (xyz)^{k_3}, (xyzxz)^{k_4}, (xzy)^{k_5}, (x^{-1}z^y)^{k_6} \}.$$

Therefore, $H$ is a homomorphic image of $G(k_1,\dots,k_6)= \langle x,y,z | B \cup R(k_1,\dots,k_6)  \rangle $, where 
$k_1,k_2,k_4, k_6 \in \{ 4,5,6,7 \}$, $k_3,k_5 \in \{4,6\}$. Note that $G(6,7,6,7,6,7)$ is a homomorphic image of $K(k_7)= \langle x,y,z | B \cup R(6,7,6,7,6,7) \cup (y \cdot (xz)^3)^{k_7}  \rangle $, where $k_7 \in \{4,6\}$. Computations show that the index in $K(6)$ of subgroup  $ \langle x, yz \rangle$, generated by an involution and an element of order $3$, is trivial. Which is not possible by Lemma 2.1. Similarly, $|K(4): \langle y,z \rangle|=1$. The group $G(6,5,6,6,6,5)$ is isomorphic to  $M_{12} \times C_3$, which is not possible. Computation show that the index in $G(6,7,6,6,6,7)$ of $ \langle x, yz \rangle$ is trivial. 

The group $H=G(6,4,6,4,6,4)$ has a normal $2$-subgroup $V$ of order $2^{11}$ and nilpotency class $2$ with $Z(V)$ of order $2^3$, and $H/V \simeq A_4$. The coset 
action of $H$ on $\langle x,yz \rangle$ provides an exact permutation representation of $H$ on $2^8$ points. Now it is straightforward to check that $H$ and its appropriate homomorphic images satisfy the conclusion of the lemma. The group $G(7,7,6,6,6,7)$ is isomorphic to $L_2(13)$, which is not possible. Computations show that the orders of other groups $G(k_1,\dots,k_6)$ divide $96$ and they are homomorphic images of $H$.

The lemma is proved. \end{proof}

\begin{wrapfigure}[120]{l}[340pt]{20mm}
	\includegraphics[width=20mm]{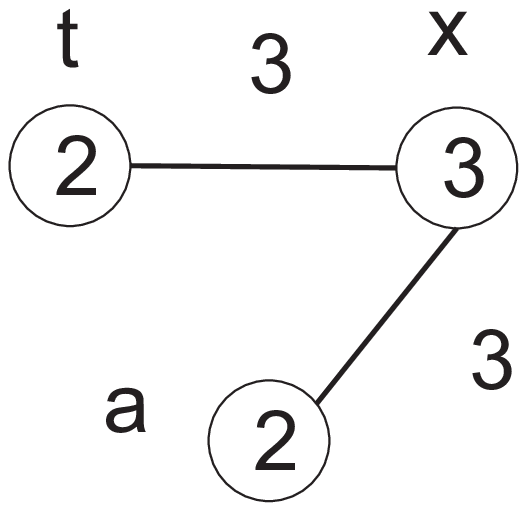}
\end{wrapfigure}

{\bf Lemma 4.5} {\it Let $t \in G$ be an involution such that $\langle t,x \rangle \simeq A_4$. Then $(at)^4=1$.}
	
\begin{proof}[Proof] Let $x_1=tx^{-1}$, $y_1=x$, $z_1=x^{-1}a$. Then $(x_1y_1)^2=(y_1z_1)^2=1$ and the statement follows from Lemma 4.4. Indeed,
$at=y_1z_1x_1y_1 \sim y_1^{-1}z_1x_1$. The lemma is proved. \end{proof}

\begin{wrapfigure}[120]{l}[340pt]{20mm}
	\includegraphics[width=20mm]{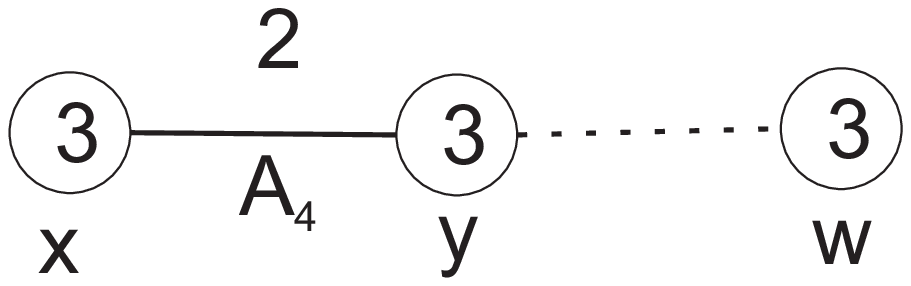}
\end{wrapfigure}

{\bf Lemma 4.6} {\it Let $x,y,w \in G$ be elements of order $3$ satisfying the following relations: $(xy)^2=[y,w]=1$. Then $\langle x,y,w \rangle$ is finite and  $((xy)(xy)^w)^4=1$.}

\begin{proof}[Proof]

Let $z=x^w$. Then  $xy \sim x^wy^w \sim x^wy \sim yz$, and hence $(yz)^2=1$. Тherefore conditions of Lemma 4.5 are true for $\langle x,y,z \rangle$. Denote by $R$ the relations from the conclusion of Lemma 4.5, and by $B$ the relations from the proof of Lemma 4.5. Then note that $\langle x,y,w \rangle$ is a homomorphic image of $G(k_1,k_2,k_3)=\langle x,y,w | R \cup B \cup \{ w^3, [y,w], (xw)^{k_1}, (x^{-1}w)^{k_2}, (xyw)^{k_3} \} \rangle$, where
$k_1,k_2,k_3 \in \{4,5,6,7 \}$. 

Note that the group $H=G(6,6,6)$ has a normal $2$-subgroup $V$ of order $2^{10}$ and nilpotency class $2$ with $Z(V)$ of order $2^4$, and $H/V \simeq C_3 \times C_3$. Computations show that orders of $H/ \langle ((xy)(xy)^w)^4 \rangle$ and $H$ coincide, meaning that the required identity holds in $H$ and
its homomorphic images. Computations show that $G(7,7,7) \simeq L_3(4)$, which is not possible, and orders of $G(k_1,k_2,k_3)$ divide $12$ in other cases. The lemma is proved. \end{proof}

{\bf Lemma 4.7} {\it Let $t\in G$ be an involution such that $(at)^3=(tx)^6=1$. Then $(tt^x)^6 \not = 1$. }

\begin{proof}[Proof] Assume the contrary. 

Note that elements $h=(t^xt)^2$ and $k=[h,x]$ of  $\langle t,x | 1=t^2=x^3=(tx)^6=(t^xt)^6 \rangle \simeq F_{216}$ satisfy the following conditions: $\langle h,x \rangle$ is an extraspecial group of period $3$ and order $27$, $k$ is its central element. Moreover, $h^t=h^{-1}$ and $k^t=k^{-1}$. 

Elements $ax^{-1}, x, k$ satisfy the conditions of Lemma 4.6, hence $(a^ka)^4=1$. Let $c=tk$. Note that $\langle a,t,c \rangle$ is a homomorphic image of $G(k_1,k_2,k_3,k_4) = \langle 1=a^2=t^2=c^2=(at)^3=(ct)^3=(a^{tc}a)^4=(ac)^{k_1}=(at \cdot c)^{k_2}=(ta^c)^{k_3}=(at \cdot a^c)^{k_4} \rangle$, where $k_1, k_3 \in \{4,5,6,7\}$, $k_2, k_4 \in \{4,6\}$. Computations show that orders of these groups divide $6$.

This way we obtain that $k=1$. It follows that $[h,x]=1$. Again using Lemma 4.6 $(a^ha)^4=1$. In an analogous way, considering $\langle a,t,h \rangle$, we obtain that $h=1$. So we may assume that $(tt^x)^2=1$.

Note that $u=tt^xtt^{x^{-1}}t$ is the central involution from $\langle t,x | 1=t^2=x^3=(tx)^6=[t,x]^2 \rangle \simeq C_2 \times A_4$. By Lemma 4.2 in $G$ the following relations hold $(au)^4=(a^xu)^4=1$. Applying Lemma 2.2 to $\langle t,a,a^x \rangle$ we obtain that the order of $ta^x$ divides $4$ or $6$. Therefore the subgroup $\langle a,t,x \rangle$ is a homomorphic image of   $G(k_1,k_2,k_3)=\langle a,t,x | 1=a^2=t^2=x^3=(at)^3=(ax)^3=(tx)^6=(t^xt)^2=(au)^4=(a^xu)^4=
(ta^x)^{k_1}=(t \cdot xa)^{k_2}=(ax \cdot t)^{k_3} \rangle$, where $k_1,k_2,k_3 \in \{4,6\}$. Computations show that the order of  $G(k_1,k_2,k_3)$ divides $12$. Which implies $t=1$. A contradiction. \end{proof}


\begin{wrapfigure}[120]{l}[340pt]{20mm}
	\includegraphics[width=20mm]{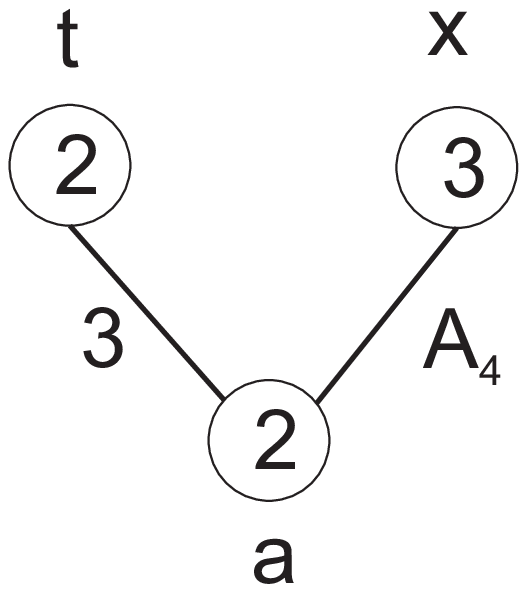}
\end{wrapfigure}

{\bf Lemma 4.8} {\it Let $t \in G$ be an involution. Then $\langle a,t \rangle \not \simeq S_3$.}

\begin{proof}[Proof] Assume the contrary.

Suppose that $bt$ and $abt$ are elements of order $6$, where $b=a^x$. Then $w=t^{bta}$ satisfies the following relations $(aw)^3=(bw)^2=(abw)^6=1$ by \cite[Лемма 5]{mamontov2013e}, and $\langle a,b,w \rangle \simeq D_{12}$.

By Lemma 2.2 changing if necessary $t$ to $w$, and $a$ to $b$ or $ab$, we may assume that $(tb)^6=(tab)^4$.

Suppose that $(tx)^4=1$. Then $\langle a,t,x \rangle$ is a homomorphic image of $G(k)=\langle a,t,x |1=a^2=t^2=x^3=(tx)^4=(ax)^3=(at)^3=(tb)^6=(tab)^4=(ax \cdot t)^k \rangle$,
where $k \in \{4,6\}$. Computations show that $G(6) \simeq V \rtimes L_2(7)$ and $G(4) \simeq A_6$, where $V$ is elementary abelian of order $2^3$. Therefore we may assume that the order of $tx$ 
is not equal to $4$.

Let $t_1=a^t$. Then $t_1$ and $t$ are conjugated with $a$ and so sit in $\Delta$ and centralize more than one involution. By Statement 2, Lemmas 4.3, 4.7, 2.1 and arguments above $(tx)^6=(t^xt)^4=1$. Similarly, $(t_1x)^6=(t_1^xt_1)^4=1$.
The involution $u=tt^xtt^{x^{-1}}t$  of the group  $\langle t,x | (tx)^6=(t^xt)^4=1 \rangle$ centralizes $x$. By Lemma 4.2 $(au)^4=1$.
Note that $\langle a,t,x \rangle$ is a homomorphic image of $G(k_1,k_2) = \langle a,t,x | 1=a^2=t^2=x^3=(ax)^3=(at)^3=(au)^4=(tb)^6=(tab)^4=(tx)^6=(t^xt)^4=(t_1x)^6=({t_1}^xt_1)^4=(ax \cdot t)^{k_1}=(xa \cdot t)^{k_2} \rangle$, where $k_1,k_2 \in \{4,6\}$. Computations show that the order of $G(k_1,k_2)$ divides $12$. Hence $t=1$ and this case is not possible.

The lemma is proved. \end{proof}

Note that we obtained a different proof of \cite[Statement 1]{mamontov2013e}.

\begin{wrapfigure}[120]{l}[340pt]{20mm}
	\includegraphics[width=20mm]{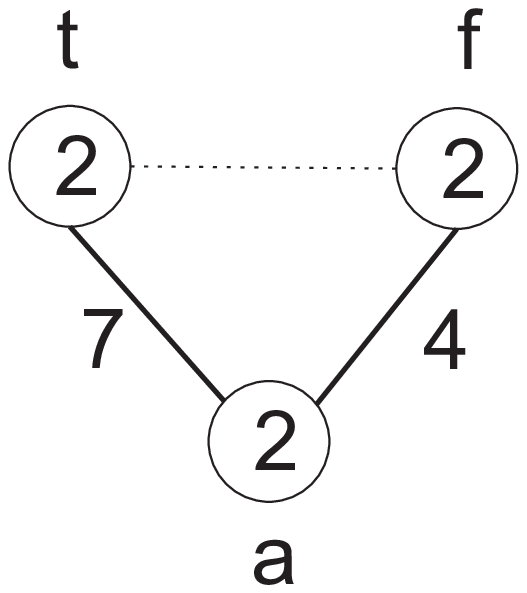}
\end{wrapfigure}

{\bf  Lemma 4.9} {\it Assume that involutions $a,t,f$ satisfy the following relations: $(at)^7=(af)^4=1$ and $[t,f]=1$. Then $(atf)^7=(t(af)^2))^4=1$ and  $\langle a,t,f \rangle$ is double frobenius $V:C_7:C_2$, where $V$ is elementary abelian of order $2^6$.}

\begin{proof}[Proof] 

To prove that the group is finite note that $\langle a,t,f \rangle$ is a homomorphic image of $G(i_1,i_2)=\langle a,t,f |1=a^2=t^2=f^2=(at)^7=(af)^4=(tf)^2=(atf)^{i_1}=(t(af)^2)^{i_2} \rangle$, where $i_1, i_2 \in \{4,5,6,7\}$. Computations show that $G(7,4)$ is finite and satisfy conclusion of the lemma, while the order of $G(i_1,i_2)$ is not greater than $14$ for other values of the parameters, which is not possible. The lemma is proved. \end{proof}

\begin{proof}[Proof of Theorem]

By Statement 2 and Lemma 4.8 we may assume that $(at)^7=1$.

First suppose that $(tx)^6=(t^xt)^2=1$. Then we may assume that $\langle t^xt ,x \rangle \simeq A_4$ and that $u=tt^xtt^{x^{-1}}t$ is an involution, which centralizes $x$. By Lemma 4.5 $(a \cdot t^xt)^4=(a^x \cdot t^xt)^4=1$. By Lemma 4.2 $(au)^4=1$. Using Lemma 4.8 we may assume that $3$ does not divide the order of $a^xt$. Consequently,  $\langle a,t,x \rangle$ is a homomorphic image of $G(k)= \langle a,t,x | 1=a^2=t^2=x^3=(ax)^3=(tx)^6=(tt^x)^2=(au)^4=(att^x)^4=(a^xtt^x)^4=(at)^{7}=(a^xt)^{k} \rangle$. Computations show that the order of  $G(k)$ divides $12$.

Now suppose that $tx \in \Gamma_6$ and $t^xt \in \Gamma_4$. Denote $y=(xt)^2$ and $v=(xt)^3$. Obviously, $[y,v]=1$. In $\langle x,t \rangle$ element $y$ of order $3$ sits inside the subgroup $\langle (t^xt)^2,y \rangle$ isomorphic to $A_4$. By Lemma 4.3 the order of $aa^y$ is not 7. Therefore, by Statement 2 and Lemma 4.8 $(aa^y)^4=1$. In $\langle x,t \rangle$ we also have $(v^xv)^2=1$. By arguments above $(av)^4=1$. Let's prove that $(aa^{xt})^4=1$.

\begin{itemize}
	\item First assume that $ay \in \Gamma_6$.
	
Let $u=aa^yayay^{-1}a$. Then $u^2=1$ and $[u,y]=1$. Having $u,v \in C_G(y)$, the order of $uv$ divides $6$.

If $\langle u,v \rangle$ has an element of order $3$, then let $w$ be an involution such that the order of $uw$ is 3. Note that if $a^g$ is an arbitrary conjugate of $a$, then the order of
$a^gv$ is even. In other case $a$ is conjugated with $v$ and hence inverts an element of order $3$, which is in contradiction with Lemma 4.8. Therefore  $\langle a,y,w \rangle$ is a homomorphic image of 
$G(i)=\langle a,y,w | 1=a^2=y^3=w^2=(ay)^6=(a^ya)^4=w^yw=(aw)^4=(uw)^3=(au^w)^4=(a^{ya}w)^4=(ayw)^{i} \rangle$.  Computations show that $G(5) \simeq S_6$ and for other $i \in \{4,5,6,7\}$  the order of $G(i)$ is not greater than $| \langle a,y \rangle | =96$. So this case is not possible.

Now suppose that $(uv)^2=1$. Then $\langle a,y,v \rangle$ is a homomorphic image of $G(i_1,i_2)=\langle a,y,v | 1=a^2=y^3=v^2=(ay)^6=[a,y]^4=[y,v]=(av)^4=(uv)^2=(ayv)^{i_1}=((av)^2y)^{i_2} \rangle$, where $i_1,i_2 \in \{4,5,6,7\}$; $u=aa^yayay^{-1}a$. Computation show that such group is always finite and contains an element of order $3$ only if $i_1=i_2=6$. In the last case $G(6,6)$ is $\{2,3 \}$-group of order $2^{18}3$, where by Lemma 4.8 we have the required equality  $(aa^{xt})^4=1$.

	\item Let now  $ay \in \Gamma_4$.
	
Then $\langle a,y,v \rangle$ is a homomorphic image of $G(i_1,i_2) = \langle a,y,v | 1=a^2=y^3=v^2=(ay)^4=[y,v]=(av)^4=(ayv)^{i_1}=((av)^2y)^{i_2} \rangle$, where $i_1,i_2 \in \{4,5,6,7\}$.
In $G(4,6)$ we have  $|aa^{yv}|=6$, which is not possible by Lemma 4.8. Further note that $G(6,4) \simeq C_2 \times S_6$, which is no possible by Statement 2. Computation show that $G(i_1,i_2)$ does not have elements of order 3 for other values of the parameters $i_1$ and $i_2$. Therefore, this case is not possible.

\end{itemize}

So further we assume that $(aa^{xt})^4=1$. Then $(a^xa^t)^4=1$ and by Lemma 4.9 we obtain that $\langle a,a^t,a^x \rangle =\langle a,t,a^x \rangle$ is the double frobenius group. In particular, $(ta^x)^4=1$ and $(taa^x)^7=(a^taa^x)^7=1$.

Note that in  $\langle t,x \rangle$ we have $(x^{-1}t)^6=1$, $\langle (x^{-1}t)^3,x \rangle \simeq C_2 \times A_4$. Therefore repeating arguments above we obtain that  $(aa^{x^{-1}t})^4=1$. Consequently,
$(a^ta^{x^{-1}})^4=(a^taa^x)^4=1$. A contradiction.

The theorem is proved. \end{proof}


\end{document}